\documentclass[11pt]{amsart}

\usepackage{amssymb,amscd,amsthm,amsxtra}
\usepackage{latexsym}
\usepackage{dsfont,latexsym, color}
\usepackage{mathrsfs}

\usepackage{color}
\definecolor{citation}{rgb}{0.2,0.58,0.2} 
\definecolor{formula}{rgb}{0.1,0.2,0.6}
\definecolor{url}{rgb}{0.3,0,0.5} 
\usepackage[colorlinks=true,linkcolor=formula,urlcolor=url,citecolor=citation]{hyperref} 
\newlength{\defbaselineskip}
\newcommand{\setlinespacing}[1]
           {\setlength{\baselineskip}{#1 \defbaselineskip}}

\newtheorem{theorem}{Theorem}[section]
\numberwithin{equation}{section}

\allowdisplaybreaks
\makeatletter
\DeclareRobustCommand*{\bfseries}{%
  \not@math@alphabet\bfseries\mathbf
  \fontseries\bfdefault\selectfont
  \boldmath
}
\makeatother

\newcommand{\Om}{\Omega}

\newcommand{\dys}{\displaystyle}

\newcommand{\dx}{\, {\rm d}x}

\def\vs{\vspace{0.9mm}}

\title[Global Compactness in fractional Sobolev spaces]
{A Global Compactness type result \\[0.2ex] for Palais-Smale sequences \\[0.2ex] in fractional Sobolev spaces}

\author[G. Palatucci]
{Giampiero Palatucci}
\address[Giampiero Palatucci]{Dipartimento di Matematica e Informatica, Universit\`a degli Studi di Parma, Campus - Parco Area delle Scienze, 53/a, 43124 Parma, Italy
\\ \qquad  \qquad \qquad \and 
SISSA, Via Bonomea, 256, 34136 Trieste, Italy}
\email{\href{mailto:giampiero.palatucci@unipr.it}{giampiero.palatucci@unipr.it}}

\author[A. Pisante]
{Adriano Pisante}
\address[Adriano Pisante]{Dipartimento di Matematica, Sapienza Universit\`a di Roma
\\ P. le Aldo Moro, 5 
\\ 00185 Roma, Italia}
\email{\href{mailto:pisante@mat.uniroma1.it}{pisante@mat.uniroma1.it}}

\begin{document}

\begin{abstract}
{We extend the global compactness result 
 by M.~Struwe (\cite{Str84}) to any fractional Sobolev spaces~$\dot{H}^s(\Omega)$, for~$0<s<N/2$ and $\Omega \subset {\mathds R}^N$ a bounded domain with smooth boundary. The proof is a simple direct consequence of the so-called {\em profile decomposition} of P.~Gerard (\cite{Ger98}).}
\end{abstract}

\subjclass[2010]{Primary  35J60, 35C20, 35B33, 49J45\vspace{1mm}}

\keywords{Profile decomposition, global compactness,  fractional Sobolev, critical Sobolev exponent.\vspace{1mm}}

\maketitle

\begin{center}
 \rule{12.7cm}{0.5pt}\\[0cm] 
{\sc {\small To appear in}\, {\it Nonlinear Anal.}}
\\[-0.25cm] \rule{12.7cm}{0.5pt}
\end{center}
\vspace{4mm}


\setlinespacing{1.02}

\section{Introduction} 

Since the seminal paper~\cite{Str84}, global compactness properties for Palais-Smale sequences in the Sobolev space~$H^1$ have become very important tools in Nonlinear Analysis which have been crucial in many existence results, e.~\!g. for ground states and blow-up solutions for nonlinear Schr\"odinger equations,  
for solutions of Yamabe-type equations both in conformal and in CR geometry, for prescribing $Q$-curvature problems, etc...
Together with the aforementioned examples concerning single equations in a scalar unknown function, more difficult systems of PDEs, often related to
 other geometric problems, share similar compactness properties for their solutions; for instance, this is the case for parametric surfaces of constant mean curvature, harmonic maps from Riemann surfaces into Riemannian manifolds, Yang-Mills connections over four-manifolds, pseudo-holomorphic curves into symplectic manifolds, planar Toda systems, etc... The involved literature is really too wide to attempt any reasonable account here. 
 \vspace{1mm}
 
 In the present note, we aim to extend the global compactness result by M.~Struwe for semilinear elliptic equation in $H^1_0$ to the case  of fractional Sobolev spaces~$\dot{H}^s$ of any real differentiability order $0<s<N/2$ by means of the so-called {\em profile decomposition} first obtained in \cite{Ger98} (see also \cite{PP13} for an alternative slightly simpler and more abstract approach).

Let $N\geq 1$ and for each $0<s<N/2$ denote by $\dot{H}^s({\mathds R}^N)$ the usual $L^2$-based homogeneous fractional Sobolev spaces\footnote{For further details on the fractional Sobolev spaces, we refer to~\cite{DPV12} and the references therein.} 
defined via Fourier transform as the completion of $C_0^\infty(\mathds{R}^N)$ under the norm 
\begin{equation}
\label{Hs0-norm}
\| u\|^2_{\dot{H}^s}=\int_{\mathds{R}^N}|\xi|^{2s}|\hat{u}(\xi)|^2 \, {\rm d}\xi \, .
\end{equation}
In view of the well known critical Sobolev embedding $\dot{H}^s \hookrightarrow L^{2^\ast}\!$, where $2^\ast=2N/(N-2s)$ is the critical Sobolev  exponent, one has equivalently, 
\[ \dot{H}^s(\mathds{R}^N)=\Big\{ u\in L^{2^{\ast}}(\mathds{R}^N) \, \,\hbox{s.\!~t.} \, \, (-\Delta )^{s/2} u \in L^2(\mathds{R}^N)\, \Big\} \, , \]
where, by definition, $((-\Delta)^{s/2}u )^{\widehat{\, \, \,}}  (\xi):=|\xi|^s \hat{u}(\xi)$.  

Let $\Om$ be a bounded domain in ${\mathds R}^N$ with smooth boundary, $N\geq1$, and define the homogeneous Sobolev space ${\dot{H}^s(\Om)}$ as the completion of $C_0^{\infty}(\Om)$ under the norm \eqref{Hs0-norm}, hence a closed subspace of $\dot{H}^s({\mathds R}^N)$.  Thus, we have a well defined fractional Laplacian $(-\Delta)^s: \dot{H}^s(\Omega) \to (\dot{H}^s(\Omega))'$ which is a bounded linear operator (isomorphism) given by $((-\Delta)^{s}u )^{\widehat{\, \, \,}}  (\xi):=|\xi|^{2s} \hat{u}(\xi)$, so that $\langle v,\, (-\Delta)^s u\rangle_{H,H^\prime}=(v,u)_H$ for any~$u,v \in H=\dot{H}^{s}(\Omega)$.

Since in the entire space $\dot{H}^s \hookrightarrow L^2_{\textrm{loc}}$ with compact embedding, we also have $\dot{H}^s(\Omega) \hookrightarrow L^2(\Omega) \hookrightarrow (\dot{H}^s(\Om))'$, both with compact embedding.
As a consequence, there is a well defined first eigenvalue $\lambda_1=\min \big\{ \| u\|^2_{\dot{H}^s} \, \, \,  u\in \dot{H}^s(\Omega) \, , \, \, \| u\|_{L^2}=1 \big\}$, with $\lambda_1=\lambda_1(\Omega)>0$ \, and also  corresponding eigenfunctions (which are positive and simple when $s\leq 1$; see, e.~\!g., \cite[Theorem 4.2]{FP14}). 
Similarly, one has an increasing sequence of positive eigenvalues (repeated with multiplicities) going to infinity $0<\lambda_1 \leq \lambda_2 \leq \ldots$ and corresponding eigenfunctions $v_1, v_2, \ldots$ giving an orthogonal base both of $L^2(\Omega)$ and of $\dot{H}^s(\Omega)$, so that $(-\Delta)^s v_k=\lambda_k v_k$ in $(\dot{H}^s(\Omega))^\prime$ for any integer $k\geq 1$, i.~\!e. $(v_k, u)_H=\lambda_k (v_k,u)_{L^2}$ for any $u\in \dot{H}^{s}(\Omega)$. Indeed it is enough to write $(u,v)_{L^2}=(Ku, v)_{\dot{H}^s}$ for some $K\in \mathcal{L}(\dot{H}^s)$ which is compact and self-adjoint and apply the spectral theorem.

\vspace{1.5mm}

For any fixed $\lambda \in \mathds{R}$, consider the following nonlinear problem
\begin{equation}\label{pb}
(-\Delta)^s u - \lambda u - |u|^{2^\ast\!-2}u=0 \ \, \text{in} \ (\dot{H}^s(\Om))' \tag{P$_\lambda$} \, ,
\end{equation}
i.~\!e. the Euler-Lagrange equation $d \mathcal{E}(u)=0$ corresponding to the differentiable functional
\begin{equation}\label{energia}
\dys
{\mathcal E}(u)= \frac{1}{2} \int_{{\mathds R}^N} |(-\Delta)^{\frac{s}{2}}u|^2\dx -\frac{\lambda}{2}\int_{\Om}|u|^{2}\dx -\frac{1}{2^\ast}\int_\Om |u|^{2^\ast}\!\dx\, . 
\end{equation}

It is worth noticing that when $\lambda<\lambda_1$, although the functional possess the Mountain Pass geometry (arguing as in \cite{Str96}, Chapter II, Section 6), the celebrated Mountain Pass lemma does not apply because the Palais-Smale condition  fails. More generally, when $\lambda_k<\lambda<\lambda_{k+1}$ the functional has a linking geometry (using the spectral decomposition above and arguing again as in \cite{Str96}, Chapter II, Section 8) but the usual minimax scheme still cannot be applied for the same reason. As it is well known when $s=1$, this is due to the presence of a limiting nonlinearity in~\eqref{energia} and it is related to the lack of compactness for the associated critical Sobolev embedding
 $\dot{H}^s \hookrightarrow L^{2^\ast}\!$,  
which is a consequence of the invariance of the $\dot{H}^s$- and $L^{2^\ast}$-norms with respect to the scaling
\begin{equation}\label{scaling}
u (\cdot) \, \rightsquigarrow \, \tilde{u}_{x_0,\eta}(\cdot) = \eta^{\frac{2s-N}{2}} u\!\left(\frac{\cdot\, \,-x_0}{\eta}\right) \, ,
\end{equation}
for arbitrarily fixed $\eta>0$ and $x_0 \in \mathds{R}^N$.

 In the seminal paper \cite{BN83} the authors circumvent this difficulty proving that, for $s=1$, a local $(PS)$-condition holds for $\lambda<\lambda_1$ small enough. Soon after a decisive breakthrough was obtained in~\cite{Str84}, still in the local case $s=1$, describing the precise mechanism responsible for the lack of the $(PS)$-condition; i.~\!e., in Author's words, proving that {\em compactness for Palais-Smale sequences holds ``apart from jumps of the topological type of admissible functions''}, a sense we will make precise below. This major advance paved the way to several extensions and to a huge number of applications, e.~\!g. in the case of problems involving biharmonic and polyharmonic operators but also in other more complicated problems (see e.~\!g.~\cite{Str96}, \cite{TF07} and the references therein).
 
In order to state precisely our main result, consider the following limiting problem
\begin{equation}\label{pb0}
(-\Delta)^s u - |u|^{2^\ast\!-2}u=0 \, \ \text{in} \ (\dot{H}^s(\Omega_0))'\,, \tag{P$_{\rm 0}$}
\end{equation}
where $\Omega_0$ either the whole~${\mathds R}^N$ or a half-space;
i.~\!e. the Euler-Lagrange equation $d\mathcal{E}^\ast(u)=0$ corresponding to the energy functional ${\mathcal E}^\ast:\dot{H}^s(\Omega_0)\to {\mathds R}$\,, 
\begin{equation}\label{energia0}
\dys
{\mathcal E}^{\ast}(u)= \frac{1}{2} \int_{{\mathds R}^N} |(-\Delta)^{\frac{s}{2}}u|^2\dx -\frac{1}{2^\ast}\int_{\Omega_0} |u|^{2^\ast}\!\dx.
\end{equation}
We have the following extension of the result in \cite{Str84}, describing Palais-Smale sequences for \eqref{energia} in the full range $0<s<{N}/{2}$,
\begin{theorem}\label{teorema}
Let ${\mathcal E}$ and ${\mathcal E}^\ast$ the functionals defined by~\eqref{energia} and~\eqref{energia0}, respectively.
Let $\{u_n\}$ be a sequence in $\dot{H}^s(\Om)$ such that
$\mathcal{E}(u_n) \leq c$ and
\begin{equation}\label{hp_forte}
 d \mathcal{E}(u_n) \, \underset{n\to\infty}\longrightarrow\,  0 \ \text{in} \ (\dot{H}^{s}(\Om))^{'}.
\end{equation}
Then, there exists a (possibly trivial) solution $u^{(0)}$ to~\eqref{pb} such that,  for a renumbered subsequence of $\{u_n\}$, we have
\[
u_n \, \underset{n\to\infty}{\rightharpoonup}  \,u^{(0)} \ \text{in} \  \dot{H}^s(\Om) \, .
\]
Moreover, either the convergence is strong, or there exist a finite set $\mathcal{J}=\{1, 2, ..., J\}$,
 nontrivial solutions~$\{u^{(j)}\}_{j \in \mathcal{J}}$ to~\eqref{pb0} either in half spaces or in the entire space, $u^{(j)} \in \dot{H}^s(\Omega_0^{(j)})$, finitely many sequences of numbers $\{\lambda_n^{(j)}\}_{j \in \mathcal{J}}\subset(0,\infty)$ converging to zero, and finitely many sequences of points~$\{x_n^{(j)}\}_{j\in \mathcal{J}} \subset \Omega$
such that,  for a renumbered subsequence of~$\{u_n\}$, we also have for any $j\in \mathcal{J}$
\[u_n^{(j)}(\cdot):= {\lambda_n^{(j)}}^{\frac{N-2s}{2}}u_n(x_n^{(j)}+\lambda_n^{(j)} \, \cdot \, ) \underset{n\to\infty}{\rightharpoonup}  \,u^{(j)}(\, \cdot \, ) \ \text{in} \ \dot{H}^s({\mathds R}^N)\]
In addition
\begin{equation}
\label{4a}
u_n(\cdot) = u^{(0)}(\cdot)+\sum_{j=1}^J  {\lambda_n^{(j)}}^{\frac{2s-N}{2}}u^{(j)} \left( \frac{\cdot-x_n^{(j)}}{\lambda_n^{(j)}}\right) +o(1) \ \text{in} \ \dot{H}^s({\mathds R}^N),
\end{equation}
and
\begin{equation}
\label{4b} 
\left| \log{\left(\frac{\lambda^{(i)}_n}{\lambda^{(j)}_n}\right)} \right| + \left|\frac{(x^{(i)}_n - x^{(j)}_n)}{\lambda^{(i)}_n}\right| \,
 \underset{n\to\infty}\longrightarrow \infty \ \ \text{for} \ i \neq j, \, i,j \in \mathcal{J}
\end{equation}
\begin{equation}\label{4c}
\dys
 \|u_n\|^2_{\dot{H}^s} \, = \, \dys \sum_{j=0}^J \|u^{(j)}\|^2_{\dot{H}^s} + o(1)\ \text{as}\  n\to \infty,
\end{equation}
\begin{equation}\label{4d}
 \mathcal{E}(u_n)  \, = \, \mathcal{E}(u^{(0)}) + \dys\sum_{j=1}^J \mathcal{E}^{\ast}(u^{(j)}) + o(1) \ \text{as}\  n\to \infty.
\end{equation}
\end{theorem}

In the local case $s=1$, the original proof of Global Compactness in~\cite{Str84} consists in a subtle analysis concerning how the Palais-Smale condition fails for the functional~${\mathcal E}$, based on rescaling arguments,  used in an iterated way to extract convergent subsequences with nontrivial limit, together with some slicing and extension procedures on the sequence of approximate solutions to \eqref{pb}. That proof is difficult to extend even to the case $s=2$ as was done in \cite{HR01} and \cite{GGS03} (see, also, the forthcoming paper~\cite{FG14} for the analog on asymptotically hyperbolic Riemannian manifolds).
The aforementioned tools seem even more cumbersome to adapt to the more general framework presented here and further difficulties appear due to the nonlocal structure of the involved fractional Sobolev spaces $\dot{H}^s$ and the corresponding equations when $s$ is not integer.

On the contrary, in the present note,  we show how to deduce the results in Theorem~\ref{teorema} in a simple way from the so-called {\em profile decomposition} for bounded sequences in~$\dot{H}^s$ spaces. The profile decomposition, as first proved by P.~G\'erard in~\cite{Ger98}, represents a far-reaching functional analytic generalization of the aforementioned result previously known only for Palais-Smale sequences and when $s$ is an  integer.
 While its essence is a somewhat general asymptotic orthogonal decomposition property, particularly transparent e.~\!g. in abstract Hilbert space setting (see~\cite[Chapter~3]{TF07}), it has been generalized to a wide range of Banach function spaces through wavelet decomposition (see \cite{BCK11} and the references therein).  In view of its extreme flexibility, it has become a common decisive tool in the study of properties of solutions of many evolution equations and related issues. 

In the following, we state such result in the form as presented in~\cite[Theorem~4]{PP13}, whose proof has been derived by the authors combining a refinement of the critical Sobolev inequality with remainder in Morrey spaces with the abstract Hilbertian approach to profile decomposition in terms of {\em dislocation spaces} developed in~\cite{TF07}.
\begin{theorem}\label{thm_dp}
Let $\{ u_n\}$ be a bounded sequence in $\dot{H}^s({\mathds R}^N)$. Then, there exist a {\rm(}at most countable{\,\rm)} set $I$, a family of profiles $\{\psi_j\}\subset \dot{H}^s({\mathds R}^N)$, a family of points $\{x^{(j)}_n \}\in {\mathds R}^N$ and a family of numbers $\{ \lambda_n^{(j)}\}\subset(0,\infty)$, such that, for a renumbered subsequence of $\{ u_n\}$, we have
\begin{equation}
\label{orthscales} 
\displaystyle
\left| \log{\left(\frac{\lambda^{(i)}_n}{\lambda^{(j)}_n}\right)} \right| + \left|\frac{(x^{(i)}_n - x^{(j)}_n)}{\lambda^{(i)}_n}\right| \,
 \underset{n\to\infty}\longrightarrow \infty \ \ \text{for} \ i \neq j, \,\, i,j \in I \, ,
\end{equation}
\begin{equation}\label{dp2}
\displaystyle
u_n(x) = \sum_{j\in I} {\lambda^{(j)}_n}^{\frac{2s-N}{2}} \psi_j\left(\frac{x-x^{(j)}_n}{\lambda_n^{(j)}}\right) + r_n(x),
\end{equation}
where\begin{equation}\label{dp2bis}
\displaystyle \,  \lim_{n\to \infty} \|r_n\|_{L^{2^\ast}} = 0,
\end{equation}\vspace{1,5mm} 
\begin{equation}\label{dp3}
\displaystyle
\hspace{-1.6cm} \text{\qquad and} 
\qquad\qquad\quad \ \ \|u_n\|^2_{\dot{H}^s} = \sum_{j\in I} \| \psi_j\|^2_{\dot{H}^s} + \|r_n\|^2_{\dot{H}^s} + o(1) 
\ \text{as} \ n \to \infty.
\end{equation}
\end{theorem}
Armed with such a strong characterization, the proof of Theorem~\ref{teorema} goes as follows. After having checked that the $(PS)$-sequence $\{u_n\}$ is bounded in $\dot{H}^s$, we analyze the corresponding profile decomposition in the sense of Theorem~\ref{thm_dp}. Thus, we prove that the weak limit solves~\eqref{pb} and all the other nontrivial profiles $\psi^{(j)}$ solve the limiting equation~\eqref{pb0} either in the entire space or in suitable half-spaces. As a consequence the profiles are in finite number. 
 Finally, we can conclude the proof by showing that the sequence of the remainders~$\{ r_n\}$ vanishes strongly in the $\dot{H}^s$-norm, which is a simple consequence of the orthogonality of scaling parameters together with an iterated use of the celebrated Brezis-Lieb lemma from~\cite{BL83}.

Finally, it is worth pointing out that when $s=1$ the case $\Omega_0$ being a half space cannot occur. This is a consequence of the standard higher regularity properties of solutions, the Pohozaev identity and the unique continuation property, which will assure that no trivial solutions of the limiting problem in such a case do exist (see also the footnote on Page~\pageref{nota} below). In the general case when $s$ is not integer, though higher regularity could be presumably deduced (see, for istance, the papers \cite{KMS14, KMS14b} for the case $0<s<1$) and remarkable Pohozaev-type identities has been recently established (see the papers by Ros-Oton and Serra~\cite{ROS14,ROS14b}), a unique continuation property is still missing e.~\!g. in the case when $s\in (0,1)$ and it is unknown to hold even for $s=2$ (see however \cite{Ped58} for the case $s=2$ in which unique continuation property is obtained under extra assumptions on the derivates of the solutions on the boundary).

\vspace{1.5mm}
In our opinion, the contribution of this note is twofold: we extend Struwe's Global Compactness result for Palais-Smale sequences to any fractional Sobolev spaces $\dot{H}^s$, hence establishing a tool which could be useful for further investigations in this nonlocal framework; moreover, since we derive it from profile decomposition, we obtain an alternative and somewhat simpler proof of such compactness result even in the classical case of integer $s=1,2$ considered in \cite{Str84} and \cite{HR01},\cite{GGS03}.

%
%

\vspace{2mm}
\section{Proof of the extended Global Compactness}

Before starting with the proof of Theorem~\ref{teorema}, we briefly fix some further notations. In the following, we denote for brevity by $\langle u , v \rangle$ and $ \, \, \dys ( u, v )$
the natural $(\dot{H}^s)^\prime-\dot{H}^s$ duality pairing and the scalar product in~$\dot{H}^s$ associated to the norm \eqref{Hs0-norm}, respectively.
Also, we recall that for any $0<s<N/2$ the following
Sobolev inequality does hold
for some positive constant~$S^{\ast}$, depending only on~$N$ and~$s$, namely
\begin{equation}\label{eq_sobolev}
\|u\|^{2^\ast}_{L^{2^\ast}} \leq S^{\ast}\|(-\Delta)^{\frac{s}{2}}u\|^{2^\ast}_{L^2} \ \ \ \forall u \in C^\infty_0(\Om),
\end{equation}
where $2^{\ast}=2N/(N-2s)$ is the critical Sobolev  exponent; 
by a density argument, the same inequality is valid on ${\dot{H}^s(\Om)}$.

Finally, in the rest of the paper we denote by $\tilde{u}_{x_0,\eta}$ the rescaling of the function~$u$ with scaling parameters $(x_0,\eta)$ as in~\eqref{scaling}; when clear from the context we do not specify explicitly the choices of $\eta$ and $x_0$ there.

%
%

\subsection{Proof of Theorem~\ref{teorema}}
For the reader's convenience, we divide the proof in a few steps.

\noindent
\\ {\it Step 1. The sequence $\{u_n\}$ is bounded in ${\dot{H}^s(\Om)}$.}
\ This is a straightforward consequence of the fact that $\{u_n\}$ is a Palais-Smale sequence for ${\mathcal E}$. Indeed
\begin{equation*} 
\left( \frac12-\frac{1}{2^\ast}\right) \int_\Omega |u|^{2^\ast}\dx 
\, = \,  {\mathcal E}(u_n) - \frac{1}{2}\langle {\rm d} {\mathcal E}(u_n),\, u_n \rangle
\,  \leq\,  |c| + o(1) \ \ \text{as} \ n \to \infty \, ,
\end{equation*}
and in turn, since $\Omega$ is bounded,   
\begin{equation}
\label{L2bound}
\int_\Omega |u_n|^2\dx \, \leq\, |\Omega|^{\frac{2s}{N}} (|c|+o(1))^{2/{2^\ast}}\ \ \text{as} \ n \to \infty  \, . 
\end{equation}
Thus
\begin{eqnarray*}
\| u_n \|^2_{\dot{H}^s} & =  & 2 {\mathcal E}(u_n) +\lambda \int_\Omega |u_n|^2 \dx+\frac{2}{2^\ast} \int_\Omega |u|^{2^\ast}\!\dx \\[1ex]
& \leq & |c| +|\lambda| \cdot |\Omega|^{\frac{2s}{N}} (|c|+o(1))^{2/{2^\ast}}+\frac{N-2s}{N}(|c|+o(1)) \ \ \text{as} \ n \to \infty \, ,
\end{eqnarray*}
which gives the claim.

\noindent
\\{\it Step 2. The weak limit (up to subsequences) $u^{(0)}$ solves~\eqref{pb}.}
\ From the previous step, passing to a subsequence we can assume that $u_n\rightharpoonup {u^{(0)}}$ weakly in ${\dot{H}^s(\Om)}$ as $n\to \infty$, thus also strongly in $L^2(\Omega)$ so that $|u_n|^{2^\ast-2}u_n \to |u^{(0)}|^{2^\ast -2} u^{(0)}$ strongly in $L^1(\Omega)$, hence weakly in $L^{(2^\ast)^\prime}(\Omega)$. Thus, by weak continuity of ${\rm d}{\mathcal E}$, the limit function~$u^{(0)}$ is a (possibly trivial) solution to ~\eqref{pb}.
If the convergence is strong we are done.

\noindent
\\ {\it Step 3. If the convergence is not strong then $\{u_n \}$ contains further profiles}.
\ We assume that $u_n\rightharpoonup {u^{(0)}}$ only weakly in ${\dot{H}^s(\Om)}$ as $n\to \infty$, thus we can apply Theorem~\ref{thm_dp} to the sequence $\{u_n\}$ and then we obtain a profile decomposition for a renumbered subsequence. We set, from now on, ${u^{(j)}}=\psi_{j}$ for a nonempty set of indices $\mathcal{J}\subset \mathds{N}$ corresponding to profiles different from $u^{(0)}$, if any. Indeed, recall that, by definition of profile, for fixed $j\in \mathcal{J}$ there exist sequences $\{x_n^{(j)}\} \subset {\mathds R}^N$ and $\lambda_n^{(j)} \subset (0,\infty)$ such that $u_n(x_n^{(j)}+\lambda_n^{(j)} \, \cdot \,) \rightharpoonup u^{(j)} (\cdot)$ in $\dot{H}^s({\mathds R}^N)$ as $n \to \infty$. Thus, if $u^{(0)}\neq 0$ then it can be regarded as profile corresponding to the trivial scaling $x_n^{(0)}\equiv 0$ and $\lambda_n^{(0)} \equiv 1$.

Here we are going to show that if the convergence is not strong then $\mathcal{J}$ is always not empty, as it always contains profiles associated to nontrivial scalings (postponing to {Step 6} the proof that $\mathcal{J}$ is also a finite set).

Arguing by contradiction, if $\{ u_n \}$ contains no profiles  \big(except possibly $u^{(0)}$\big) then $u_n \to u^{(0)}$ strongly in $L^{2^\ast}$\,(see e.~\!g. \cite[Proposition 1]{PP13}). On the other hand
\begin{eqnarray*}
\| u_n-u^{(0)}\|^2_{\dot{H}^s} & = & ( u_n, u_n-u^{(0)}) -(u^{(0)},u_n-u^{(0)}) 
\\[0.5ex] 
& = & \langle d\mathcal{E}(u_n), u_n-u^{(0)}\rangle 
- \langle d\mathcal{E}(u^{(0)}), u_n-u^{(0)}\rangle+\lambda \| u_n-u^{(0)}\|_{L^2}^2 \\[0.5ex]
 && \dys +\int_\Omega \left( |u_n|^{2^\ast-2}u_n-|u^{(0)}|^{2^\ast-2}u^{(0)}\right) \left( u_n-u^{(0)}\right) \dx
 \ = \ o(1)
\end{eqnarray*}
as $n \to \infty$ in view of \eqref{hp_forte}, {Step 2} and the strong convergence $u_n \to u^{(0)}$ in $L^{2^\ast}$. Thus the contradiction shows that $\mathcal{J}\neq \emptyset$, i.~\!e. $\{ u_n \}$ contains further profiles.

\noindent
\\ {\it Step 4. The scaling parameters from Theorem \ref{thm_dp} satisfy $\ \lambda_n^{(j)}\to 0$ and \,$\text{\rm dist}\big(x_n^{(j)}, \,\Omega\big)=\mathcal{O}(\lambda_n^{(j)})$ as $n \to \infty$ for any $j\geq 1$.} 
\ We argue by contradiction and we first assume that, up to a subsequence, $\lambda_n^{(j)} \to \infty$.
Since the convergence defining the profiles is actually strong in $L^2_{\textrm{loc}}$ and $u^{(j)} \neq 0$, we have
\begin{eqnarray*}
0 & < & \int_{{\mathds R}^N}|u^{(j)}|^2\dx\\
& \leq & \liminf_n \int_{{\mathds R}^N} {\lambda_n^{(j)}}^{N-2s}|u_n(x_n+ \lambda_n^{(j)} \, \cdot \,)|^2 \dx 
\ = \ \liminf_n {\lambda_n^{(j)}}^{-2s} \int_\Omega |u_n|^2 \dx
\end{eqnarray*}
which, in view of \eqref{L2bound},  yields a contradiction. Thus $\lambda_n^{(j)}\leq c^{(j)}$ for each $n\geq 1$. 
Now assume that, up to subsequences, $\lambda_n^{(j)} \geq c_j>0$   for each $n\geq1$.  Then, up to subsequences, $|x_n| \to \infty$, since otherwise we would clearly have $u^{(0)} \neq 0$, which is however impossible because the two profiles $u^{(0)}$ and $u^{(j)}$, $j\geq 1$, must obtained along distinct sequence of scaling parameters  satisfying~\eqref{orthscales}. We conclude that $|x_n|\to \infty$, but then $u_n(x_n^{(j)}+\lambda_n^{(j)} \, \cdot \,) \rightharpoonup 0$ since for fixed $n\geq 1$ each function is supported in ${\lambda_n^{(j)}}^{-1} \left(\Omega-x_n \right)$, so they are eventually zero on each compact set in $\mathds{R}^N$ since $\Omega$ is bounded and $\lambda_n^{(j)}\geq c_j>0$. Since $u^{(j)}\neq 0$ we have a contradiction and therefore $\lambda_n^{(j)} \to 0$ as $n \to \infty$.
Finally, assuming \,$\text{\rm dist}\big(x_n^{(j)},\, \Omega\big)/ \lambda_n^{(j)} \to \infty$ possibly along a subsequence, as $n \to \infty$ would yield \,$\text{\rm dist}\Big( 0,\, {\lambda_n^{(j)}}^{-1}\!\big(\Omega-x_n \big) \Big) \to \infty$ and still $u_n(x_n^{(j)}+\lambda_n^{(j)} \, \cdot \,) \rightharpoonup 0$ arguing as above, which is impossible.

\noindent
\\ {\it Step 5. The profiles $u^{(j)}$ solve~\eqref{pb0} either in a half space or in the entire space\footnote{\ In the full range $0<s<N/2$ \label{nota} the limiting domain $\Omega_0$ can be either the whole ${\mathds R}^N$ or a half space. Indeed, as mentioned in the introduction, when $s=1$ we can exclude the existence of nontrivial solutions to the limiting problems in the half space by Pohozaev identity and unique continuation. Similarly, as shown in~\cite[Theorem 1.5]{FW12}, this is still the case when dealing with {\it nonnegative solutions} in the  
 range $s\in (0,1)$. However, even when $s=2$ this possibility cannot be apriori excluded.} and one may assume $\{ x_n \} \subset \Omega$.}
 
 \ This claim follows by the assumption in~\eqref{hp_forte} together with the invariance of the $\dot{H}^s$- and $L^{2^\ast}$\!-norm with respect to the scaling in~\eqref{scaling}. According to {Step 4} we just need to distinguish two cases, depending whether $\text{\rm dist}\big(x_n^{(j)}, \,\partial  \Omega\big)=\mathcal{O}(\lambda_n^{(j)})$ as $n\to \infty$ (Case~I) or, up to a subsequence, $\text{\rm dist}\big(x_n^{(j)},\, \partial  \Omega\big)/\lambda_n^{(j)}\to \infty$ (Case~II). Note that the latter happens only when $\{ x_n^{(j)}\} \subset \Omega$ stay inside the domain or approaches the boundary ``slower than $\lambda_n^{(j)}$''.  

Suppose that $u^{(j)}_n( \cdot):= u_n(x_n^{(j)}+\lambda_n^{(j)} \, \cdot \,) \rightharpoonup u^{(j)} (\cdot)$ in $\dot{H}^s({\mathds R}^N)$ as $n \to \infty$ where, as in {Step 4}, $u^{(j)}_n \in \dot{H}^s\big({\lambda_n^{(j)}}^{-1} \left(\Omega-x_n \right)\big)$. Note that because of the condition on the scaling parameters, arguing as in the appendix of \cite{GGS03}, ${\lambda_n^{(j)}}^{-1} \left(\Omega-x_n \right) \to \Omega_0^{(j)}$ as $n\to \infty$, where $\Omega_0^{(j)}$ is either an open half space because of the smoothness of $\partial \Omega$ (Case~I) or the entire space (Case~II), in the sense that for any compact set $K \subset \Omega_0^{(j)}$ we have $K\subset {\lambda_n^{(j)}}^{-1} \left(\Omega-x_n \right)$ for all $n$ large enough. Since for the complements ${\lambda_n^{(j)}}^{-1} \left(\mathcal{C}\Omega-x_n \right) \to \mathcal{C}\Omega_0^{(j)}$ as $n\to \infty$, we also have $u^{(j)}=0$ a.~\!e. on $\mathcal{C}\overline{\Omega_0^{(j)}}$, hence $u^{(j)}\in \dot{H}^s(\Omega_0^{(j)})$.

Fix $\phi \in C^\infty_0(\Omega_0^{(j)})$. By using the invariance in~\eqref{scaling}, with $\eta=\lambda_n^{(j)}$ and $x_{\rm o}=x^{(j)}_n$ there, together with the previous convergence we have $\tilde{\phi} \in C^\infty_0(\Omega)$ and
we can write\begin{equation}\label{eq_test}
 \langle {\rm }{\rm d}{\mathcal E}(u_n), \, {\lambda_n^{(j)}}^{\frac{2s-N}{2}}\phi\!\left(\frac{\cdot \, \, -x_n}{\lambda^{(j)}_n}\right)\rangle  = \, \langle {\rm d}{\mathcal E}^\ast(u^{(j)}_n) , \, \phi \rangle - \lambda  \left( \lambda^{(j)}_n \right)^{2s} \int_{{\mathds R}^N}  u^{(j)}_n  \, \phi \dx.
\end{equation}
In view  of~\eqref{hp_forte} the left hand-side in~\eqref{eq_test} goes to zero as $n \to\infty$, since $\| \tilde{\phi}\|_{\dot{H}^s}=\| \phi\|_{\dot{H}^s}$. On the other hand, since $u^{(j)}_n \to u^{(j)}$ in $L^2_{\textrm{loc}}$ as $n\to \infty$ and $\phi$ is smooth and compactly supported, the last term in \eqref{eq_test} vanish as $n\to \infty$ because $\lambda_n^{(j)} \to 0$ from the previous step. Thus, by weak continuity of ${\rm d}{\mathcal E}^\ast$ we infer $\langle{\rm d}{\mathcal E}^\ast (u^{(j)}),\, \phi\rangle=0$  for each $\phi \in C^\infty_0(\Omega_0^{(j)})$ and arguing by density the claim follows.  Finally, we claim that even in Case I we may always assume $\{x_n^{(j)}\} \subset \Omega$. Indeed, if we fix $\bar{x}^{(j)} \in \Omega_0^{(j)}$ and we set $\bar{x}_n^{(j)}=x^{(j)}_n+\lambda_n^{(j)} \bar{x}^{(j)}$, then $\{ \bar{x}_n^{(j)}\} \subset \Omega$ for $n$ large enough and $\bar{u}^{(j)}_n( \cdot):= u_n(\bar{x}_n^{(j)}+\lambda_n^{(j)} \, \cdot \,)= u_n(x_n^{(j)}+\lambda_n^{(j)} \, ( \, \cdot \, +\bar{x}^{(j)})) \rightharpoonup u^{(j)} (\, \cdot \, +\bar{x}^{(j)})$ in $\dot{H}^s({\mathds R}^N)$ as $n \to \infty$.

Thus, the claim follows taking $\bar{x}_n^{(j)}$ as new traslation parameters and $\bar{u}^{(j)}(\cdot ):= u^{(j)}(\,  \cdot \, + \bar{x}^{(j)})$ as new profiles.

\noindent
\\ {\it Step 6. The profiles $u^{(j)}$ are in finite number.}
\ 
This is a consequence of~\eqref{dp3} in Theorem~\ref{thm_dp}. It will just suffice to show that the $\dot{H}^s$-energy of the nontrivial profiles~$u^{(j)}$, $j\geq 1$, is uniformly bounded from below. For this, we can use the previous step by testing~\eqref{pb0} with $\phi=u^{(j)}$. From the Sobolev inequality~\eqref{eq_sobolev} we get
$$
\|u^{(j)}\|^2_{\dot{H}^s} \, = \, \|u^{(j)}\|^{2^\ast}_{L^{2^\ast}} \, \leq \, S^{\ast}\left( \|u^{(j)}\|^2_{\dot{H}^s}\right)^{\frac{2^\ast}{2}}.
$$
As a consequence $\|u^{(j)}\|_{\dot{H}^s} \, \geq \, {S^\ast}^{\frac{2}{2-2^\ast}}$ and the sum in~\eqref{dp3} must be finite.

\noindent
\\ {\it Step 7. The sequence of remainders $\{r_n\}$ converges strongly to 0 in ${\dot{H}^s({\mathds R}^N)}$.}
\ In view of the asymptotic orthogonality of the scaled profiles corresponding to \eqref{orthscales} and the invariance of the $\dot{H}^s$-norm we have
\begin{eqnarray}\label{eq_abc}
\dys
\left\|r_n \right\|^2_{\dot{H}^s} & = & \left\|u_n - u^{(0)}\right\|^2_{\dot{H}^s} + \sum_{j=1}^J\left\|{\lambda^{(j)}_n}^{\frac{2s-N}{2}} u^{(j)}\left(\frac{\cdot-x^{(j)}_n}{\lambda_n^{(j)}}\right)\right\|^2_{\dot{H}^s} \nonumber \\
& & -2 \left( u_n-u^{(0)}, \, \sum_{j=1}^J {\lambda^{(j)}_n}^{\frac{2s-N}{2}} u^{(j)}\left(\frac{\cdot-x^{(j)}_n}{\lambda_n^{(j)}}\right) \right) +o(1) \nonumber \\[1ex]
& = &
\left\|u_n - u^{(0)}\right\|^2_{\dot{H}^s} + \sum_{j=1}^J \left\| u^{(j)} \right\|^2_{\dot{H}^s}+o(1)
\nonumber\\
&& -2\sum_{j=1}^J \left( u_n^{(j)}, u^{(j)} \right)  +
2\sum_{j=1}^J \left( u^{(0)},\,  {\lambda^{(j)}_n}^{\frac{2s-N}{2}} u^{(j)}\left(\frac{\cdot-x^{(j)}_n}{\lambda_n^{(j)}} \right) \right)
\nonumber\\[1ex]
& = & \left\|u_n - u^{(0)}\right\|^2_{\dot{H}^s} - \sum_{j=1}^J\int_{{\mathds R}^N}|\tilde{u}^{(j)}|^{2^\ast}\!\dx +o(1) \,.  
\end{eqnarray}
where in the last equality we also used the fact that the profiles solve~\eqref{pb0} because of {Step 5}, hence $\| u^{(j)}\|_{\dot{H}^s}^2=\int |u^{(j)}|^{2^\ast}\!\dx$ for each $j \in \mathcal{J}$.
\vspace{1.5mm}
On the other hand, since $u_n \to u^{(0)}$ strongly in $L^2$ as $n\to \infty$, arguing as in {Step 3} we have 
\[
\| u_n-u^{(0)}\|^2_{\dot{H}^s}  =  \int_\Omega \left( |u_n|^{2^\ast-2}u_n-|u^{(0)}|^{2^\ast-2}u^{(0)}\right) \left( u_n-u^{(0)}\right) \dx
 + o(1).
\]
As in {Step 2}, we have also $u_n\rightharpoonup {u^{(0)}}$ weakly in $L^{2^\ast}(\Omega)$ and  $|u_n|^{2^\ast-2}u_n \rightharpoonup |u^{(0)}|^{2^\ast -2} u^{(0)}$ weakly in $L^{(2^\ast)^\prime}(\Omega)$ as $n\to \infty$, therefore the previous equality becomes

\[
\| u_n-u^{(0)}\|^2_{\dot{H}^s}  =\int_\Omega |u_n|^{2^\ast} \dx- \int_\Omega |u^{(0)}|^{2^\ast} \dx
 + o(1) \, ,
\]
as $n \to \infty$, which combined with \eqref{eq_abc} finally gives

\begin{equation}
\label{resto}
\lim_{n\to \infty } \left\|r_n \right\|^2_{\dot{H}^s}=  \lim_{n\to \infty} \int_\Omega |u_n|^{2^\ast} \dx- \int_\Omega |u^{(0)}|^{2^\ast} \dx- \sum_{j=1}^J\int_{{\mathds R}^N}|\tilde{u}^{(j)}|^{2^\ast}\!\dx
 \, .
\end{equation}

Finally, it remains to show that
\begin{equation}
\label{brezislieb}
\int_{\Om}|u_n|^{2^\ast}\!\dx \,\underset{n\to\infty}\longrightarrow  \,
\int_{\Om}|u^{(0)}|^{2^\ast}\!\dx +
\sum_{j=1}^J\int_{{\mathds R}^N} | u^{(j)}|^{2^\ast}\!\dx.
\end{equation}
By~\eqref{dp2}, we can write $u_n = u^{(0)}+\sum_j{\lambda_{n}^{(j)}}^{2s-N/2}u^{(j)}(\, \cdot -x_n^{(j)}/\lambda_n^{(j)}) + r_n$, and if we use iteratively the Brezis-Lieb lemma and the invariance of the $L^{2^\ast}$\!-norm together with~\eqref{dp2bis} the conclusion follows (alternatively, see \cite[Equation~(1.11)]{Ger98}, which gives the same property for the general case of countably many profiles).

\noindent
\\ {\it Step 8. Completion of the proof.} \ 
Clearly \eqref{4b} comes from \eqref{orthscales} and \eqref{4a}, \eqref{4c} follow from \eqref{dp3} in view of {Step~7}. It remains to show the validity of~\eqref{4d}. Recall that $u_n \to u^{(0)}$ in $L^2(\Omega)$ and  $u^{(0)}$ solves $(P_\lambda)$ because of {Step~2}. Recall also that $\| u^{(j)}\|_{\dot{H}^s}^2=\int |u^{(j)}|^{2^\ast}\!\dx$ because of {Step~5}. Since 
\begin{equation}\label{conv}
{\mathcal E}(u_n) = \frac{1}{2}\int_{{\mathds R}^N} |(-\Delta)^{\frac{s}{2}} u_n|^2\dx -\frac{\lambda}{2}\int_{\Om}|u_n|^{2^\ast}\dx - \frac{1}{2^\ast}\int_{\Om}|u_n|^{2^\ast}\!\dx \, ,
\end{equation}
combining \eqref{4c} and \eqref{brezislieb} the conclusion follows easily as $n\to \infty$.
\hfill \mbox{} \qed

\vspace{1mm}
\noindent
\\ {\bf Acknowledgements.}  The authors would
 like to thank Marco Squassina for having pointed out to them the paper~\cite{FW12}. The first author has been supported by PRIN 2010-2011 ``Calcolo delle Variazioni''. The first author is member of Gruppo Nazionale per l'Analisi Matematica, la Probabilit\`a e le loro Applicazioni (GNAMPA) of Istituto Nazionale di Alta Matematica ``F.~Severi'' (INdAM).

\vspace{2mm}

\end{document}